\documentclass[11pt,dvips]{amsart}  
\usepackage{amssymb,amsmath,amsfonts}  
\usepackage[dvips]{graphicx,color}
\usepackage{cite}
\allowdisplaybreaks[1]
\definecolor{pgray}{gray}{0.8}

\newtheorem{theorem}{Theorem}[section]
\newtheorem{cor}[theorem]{Corollary}
\newtheorem{proposition}[theorem]{Proposition}
\newtheorem{lemma}[theorem]{Lemma}

\newtheorem{Def}[theorem]{Definition}

\numberwithin{equation}{section}

\title{\mbox{}}

\begin{document}
\begin{center}
{\bf \LARGE{
	Time Evolution of the Navier-Stokes Flow in Far-field
}}\\
\vspace{5mm}
Masakazu Yamamoto (Graduate School of Science and Technology, Niigata University)
\end{center}
\maketitle
\vspace{-15mm}
%
\begin{abstract}
Asymptotic expansion in far-field for the incompressive Navier-Stokes flow are established.
Under moment conditions on the initial vorticity, technique of renormalization together with Biot-Savard law derives an asymptotic expansion for the velocity with high-order.
Especially scalings and large-time behaviors of the expansions are clarified.
By employing them, time evolution of velocity in far-field is drawn.
As an appendix, asymptotic behavior of solutions as time variable tends to infinity is given.
\end{abstract}

\section{Introduction}
We consider spatial decay of solutions to the incompressible Navier-Stokes equations in whole space.
In several preceding works, decay-rate of solutions as $t \to + \infty$ and as $|x| \to + \infty$ are established by deriving asymptotic expansions.
Those expansions require fast-decay for the initial data and solutions.
However, even if the initial data decays fast, the velocity decays slowly as $|x| \to +\infty$.
This structure of solutions disturbs to derive asymptotic expansions with high-order.
By studying the related vortex equation and employing Biot-Savard law, we avoid this difficulty.
For simplicity, we treat only two dimensional case.
Here we study the following initial-value problem:
\begin{equation}\label{NS}
\left\{
\begin{array}{lr}
	\partial_t u + u\cdot\nabla u = \Delta u - \nabla p,
	&
	t>0,~ x \in \mathbb{R}^2,\\
	\nabla\cdot u = 0,
	&
	t>0,~ x \in \mathbb{R}^2,\\
	u(0,x) = a (x),
	&
	x \in \mathbb{R}^2,
\end{array}
\right.
\end{equation}
where $u = (u_1,u_2)$ and $p$ denote unknown velocity and pressure, respectively.
The solenoidal condition $\nabla \cdot a = 0$ also is assumed for the given initial velocity $a$.
Well-posedness, smoothness and global existence on time of solutions are very important for this problem.
Those questions are solved in several studies (for example, see \cite{FrwgKznShr08,FjtKt,Gg-Mykw,Gg-Mykw-Osd,Kzn89JDE,KznOgwTnuch02,Wisslr} and references therein).
In this paper, we treat a smooth and global solution $u$ which satisfies that
\begin{equation}\label{decay}
	\| u(t) \|_{L^1 (\mathbb{R}^2)} \le C t^{- \frac12}
	\quad\text{and}\quad	
	\| u(t) \|_{L^q (\mathbb{R}^2)} \le C (1+t)^{-\gamma_q - \frac12}
\end{equation}
for $1 < q \le \infty$ and any $t > 0$, where $\gamma_q = 1 - \frac1q$ is the decay-rate of two dimensional gaussian in $L^q (\mathbb{R}^2)$.
For why the case $q = 1$ needs the special treatment, see the sentences under proof of Lemma \ref{lem-asymplinsp} in Section 2.
Those estimates are confirmed under several frameworks by applying the solenoidal condition (cf.\cite{Brndls,Kt,Lry,Mykw98FE,Schnbk86,Schnbk91,Wgnr87,Wgnr94}), and give the upper bound of decay-rate of velocity as $t \to +\infty$.
The lower bounds are established by an asymptotic expansion.
It is well-known that fast decay of $u$ is required to introduce the asymptotic expansion.
For the heat equation, spatial decay of solutions are inherited from an initial data.
Whereas for \eqref{NS}, decay of $u$ as $|x| \to +\infty$ is not controlled by $a$.
More precisely, even if $a \in C_0 (\mathbb{R}^2)$, then
\begin{equation}\label{sp-decay}
	u(t,x) = O (|x|^{-3})
\end{equation}
as $|x| \to +\infty$ for any fixed $t > 0$ (cf.\cite{AmrchGrltSchnbk,Brndls-Vgnrn,Mykw02FE}).
Since the asymptotic expansions introduced in the preceding works need fast decay of $u$ as $|x| \to +\infty$, the polynomial decay \eqref{sp-decay} is troublesome.
Similar problem is appearing in several equations which contain nonlocal operators.
For example, far-field asymptotics of solutions to a semi-linear anomalous diffusion equation are studied in \cite{Brndls-Krch,Ishg-Kwkm-Mchhs,Y}.
In those works, spatial decay of solutions are clarified by using the idea of spatial renormalization.
Now we are interested to far-field asymptotics of the velocity $u$.
The spatial renormalization may solve it.
However, if we choose this tool, then estimates should be complicated.
In this paper, we employ the related vorticity instead of the spatial renormalization.
The vorticity is given by $\omega = \nabla \times u = \partial_1 u_2 - \partial_2 u_1$ and fulfills that
\begin{equation}\label{vort}
\left\{\begin{array}{lr}
	\partial_t \omega + u \cdot \nabla\omega = \Delta \omega,	&	t > 0,~ x \in \mathbb{R}^2,\\
	\Delta u = \nabla^\bot \omega,						&	t > 0,~ x \in \mathbb{R}^2,\\
	\omega (0,x) = \omega_0 (x),						&	x \in \mathbb{R}^2,
\end{array}\right.
\end{equation}
where $\nabla^\bot = (-\partial_2,\partial_1)$ and $\omega_0 = \nabla \times a$.
From the definition and the solenoidal condition, it is natural that $\int_{\mathbb{R}^2} x^\alpha \omega_0 (x) dx = 0$ for $|\alpha| \le 1$.
Note that the vorticity is a scholar-field in two dimensional case and restores the velocity through Biot-Savard law:
\begin{equation}\label{BSlaw}
	u = -\nabla^\bot (-\Delta)^{-1} \omega.
\end{equation}
We emphasize that decay of the vorticity as $|x| \to +\infty$ is contralled by the initial vorticity $\omega_0$.
Indeed it is satisfied that
\begin{equation}\label{decay-vort}
	\| \omega (t) \|_{L^q (\mathbb{R}^2)} \le C (1+t)^{-\gamma_q-1}
	\quad\text{and}\quad
	\| |x|^k \omega (t) \|_{L^q (\mathbb{R}^2)}
	\le C t^{-\gamma_q} (1+t)^{-1+\frac{k}2}
\end{equation}
for $1 \le q \le \infty$ and some $k \in \mathbb{Z}_+$ (see for example \cite{Kkvc,Kkvc-Trrs}).
Those estimates originally are derived from the energy methods.
The first and second inequalities seem to discontinuous for $k$.
More precisely, the second inequality with $k = 0$ has the extra singularity as $t \to +0$.
In \eqref{decay-vort}, we set the simple estimate for $|x|^k \omega$ since a natural singularity is little complicated (see the sentences before References).
Reader may confirm them by the $L^p$-$L^q$ estimates which are in proof of Proposition \ref{prop-AsympVort} in Section 2 with Gronwall type technique.
Those estimates suggest that, on $\omega = \partial_1 u_2 - \partial_2 u_1,~ \partial_1 u_2$ and $\partial_2 u_1$ are canceled in far-field.
This structure of the vortex will play same role as the spatial renormalization in our main results.
From \eqref{decay-vort} we can define an asymptotic expansion of $\omega$ in far-field with arbitrary high-order.
Since the velocity is connected to the vorticity through Biot-Savard law, we expect that an asymptotic expansion with high-order is determined also for $u$.
Those idea firstly are established by Kukavica and Reis \cite{Kkvc-Ris}, and they showed the following estimate:

For $2 \le q \le \infty,~ m \ge 2$ and $0 \le \mu < m + 2(1-\frac1q)$,
\[
	\biggl\| |x|^\mu \biggl( u(t)
	+ \sum_{2 \le |\alpha| \le m} \frac{\nabla^\alpha \nabla^\bot (-\Delta)^{-1} G(t)}{\alpha!} \int_{\mathbb{R}^2} (-y)^\alpha \omega (t,y) dy \biggr) \biggr\|_{L^q (\mathbb{R}^2)}
	=
	O ( t^{-\gamma_q - \frac12 + \frac\mu{2}} )
\]
as $t \to +\infty$.
Here $x^\alpha \omega_0 \in L^1 (\mathbb{R}^2)$ are assumed.

This estimate describes the asymptotic expansion of $u$ in far-field with arbitrary high-order since
$
	|x|^\mu \nabla^\alpha \nabla^\bot (-\Delta)^{-1} G(t) \not\in L^q (\mathbb{R}^2)
$
for large $\mu$.
Namely this estimate suggests that $u(t)$ and the summation of $\nabla^\alpha \nabla^\bot (-\Delta)^{-1} G(t)$ are cancelled in far-field.
Particularly, since $\nabla^\alpha \nabla^\bot (-\Delta)^{-1} G(t) = O (|x|^{-3})$ as $|x| \to +\infty$ when $|\alpha|=2$ (see Lemma \ref{RzG}), we see from this estimate that \eqref{sp-decay} is essential.
On the other hand, large-time behavior of $u$ is covered yet.
Especially time evolutions of the coefficients $\int_{\mathbb{R}^2} (-y)^\alpha \omega (t,y) dy$ are not clarified.
The assertion of our main theorem will solve them.
By applying Duhamel principle and the solenoidal condition to \eqref{vort}, we obtain that
\begin{equation}\label{mcurl}
	\omega (t) = G(t) * \omega_0 - \int_0^t \nabla G (t-s) * (\omega u) (s) ds,
\end{equation}
where $*$ denotes the convolution in space.
Since $\int_{\mathbb{R}^2} x^\alpha \omega_0 (x) dx = 0$ for $|\alpha| \le 1$ is assumed and $\int_{\mathbb{R}^2} (\omega u) (t,x) dx$ $= 0$ holds, $\omega$ decays fast as $t \to + \infty$ as \eqref{decay-vort}.
Also in far-field, $\omega$ decays fast (see Proposition \ref{prop-AsympVort}).
Biot-Savard law yields from \eqref{mcurl} that
\begin{equation}\label{msu}
	u(t) = -\nabla^\bot (-\Delta)^{-1} G(t) * \omega_0 + \int_0^t R^\bot R G (t-s) * (\omega u) (s) ds.
\end{equation}
Here $R^\bot R = \nabla^\bot (-\Delta)^{-1} \nabla$ is the tensor operator combined from the Riesz transforms.
More precisely, we define $R^\bot R$ as following.
%
\begin{Def}\label{def-Rz2}
For $f = (f_1,f_2) : \mathbb{R}^2 \to \mathbb{R}^2$ and $g : \mathbb{R}^2 \to \mathbb{R}$, we denote that
\[
	R^\bot R g * f
	=
	R^\bot (Rg * f)
	=
	\left( - R_2 R_1 g * f_1 - R_2^2 g * f_2,\, R_1^2 g * f_1 + R_1 R_2 g * f_2 \right),
\]
though, for a constant $b = (b_1,b_2) \in \mathbb{R}^2$,
\[
	R^\bot Rg b
	=
	R^\bot (Rg \cdot b)
	=
	\left( -(R_2 R_1 g)b_1 - (R_2^2 g) b_2,\, (R_1^2 g)b_1 + (R_1 R_2 g) b_2 \right),
\]
where $R_1 = \partial_1 (-\Delta)^{-1/2}$ and $R_2 = \partial_2 (-\Delta)^{-1/2}$ are the Riesz transforms.
\end{Def}
If we choose the vertical basis, then $R^\bot R$ is arranged as
\[
	R^\bot R
	=
	\begin{pmatrix}
		-R_2 R_1	&	-R_2^2\\
		R_1^2	&	R_1 R_2
	\end{pmatrix}.
\]
Therefore we see Hausdorff-Young inequality $\| R^\bot R g * f \|_{L^q (\mathbb{R}^2)} \le \| R^\bot R g \|_{L^p (\mathbb{R}^2)} \| f \|_{L^r (\mathbb{R}^2)}$ for $1 + \frac1q = \frac1p + \frac1r$ with $p > 1$.

The integral equation \eqref{msu} is equivalent to the usual form
\begin{equation}\label{HFK}
	u(t) = G(t) * a - \int_0^t \nabla G(t-s) * P (u \otimes u) (s) ds,
\end{equation}
where $a = - \nabla^\bot (-\Delta)^{-1} \omega_0$ and $P$ is Helmholtz-Fujita-Kato projection.
Both \eqref{msu} and \eqref{HFK} contain the Riesz transforms in their nonlinear terms.
Due to effects of them, $u$ decays slowly in far-field as \eqref{sp-decay}.
Throughout this paper, we adopt \eqref{msu}.
The asymptotic expansion of $u$ with low-order is given by summation of
\begin{equation}\label{Um}
	U_m (t)
	=
	-\sum_{|\alpha| = m+1} \frac{\nabla^\alpha \nabla^\bot (-\Delta)^{-1} G(t)}{\alpha!}
	\int_{\mathbb{R}^2} (-y)^\alpha \omega_0 (y) dy
\end{equation}
and
\begin{equation}\label{Uminf}
	U_m^\infty (t)
	=
	\sum_{|\beta|=m} \frac{\nabla^\beta R^\bot RG (t)}{\beta!}
	\int_0^\infty \int_{\mathbb{R}^2} (-y)^\beta (\omega u) (s,y) dyds
\end{equation}
for $m = 1$ and $2$.
Indeed, the following estimates are known.
\begin{proposition}\label{prop-lowt}
	Let $\omega_0 \in L^1 (\mathbb{R}^2) \cap L^\infty (\mathbb{R}^2),~ |x|^3 \omega_0 \in L^1 (\mathbb{R}^2)$ and $\int_{\mathbb{R}^2} x^\alpha \omega_0 (x) dx = 0$ for $|\alpha| \le 1$.
	Assume that the solutions $u$ of \eqref{NS} for $a = - \nabla^\bot (-\Delta)^{-1} \omega_0$ and $\omega$ of \eqref{vort} fulfill \eqref{decay} and \eqref{decay-vort} for $k = 3$, respectively.
	Then
	\[
		\Bigl\| u(t) - \sum_{m=1}^2 (U_m+U_m^\infty) (t) \Bigr\|_{L^q (\mathbb{R}^2)}
		= o ( t^{-\gamma_q-1} )
	\]
	as $t \to +\infty$ for $1 \le q \le \infty$.
	In addition, if $|x|^4 \omega_0 \in L^1 (\mathbb{R}^2)$, then
	\begin{equation}\label{lwlg}
		\Bigl\| u(t) - \sum_{m=1}^2 (U_m+U_m^\infty) (t) \Bigr\|_{L^q (\mathbb{R}^2)}
		= O ( t^{-\gamma_q-\frac32} \log t )
	\end{equation}
	as $t \to +\infty$ for $1 \le q \le \infty$.
\end{proposition}
Those estimates are developed by Carpio \cite{Crpo}, and Fujigaki and Miyakawa \cite{Fjgk-Mykw} essentially.
We find $u$ and $\omega$ which satisfy \eqref{decay} and \eqref{decay-vort}, respectively, at least if $\omega_0$ is sufficiently smooth or small (see \cite{Crpo2,Gg-Mykw-Osd,Kt,Kkvc,Kkvc-Trrs} and also Lemmas \ref{lem-asymplint} and \ref{lem-asymplinsp} in Section 2).
Since the scaling properties $\lambda^{2+m} (U_m, U_m^\infty) (\lambda^2 t, \lambda x) = (U_m, U_m^\infty) (t,x)$ for $\lambda > 0$ are fulfilled, Proposition \ref{prop-lowt} yields the large-time behavior of $u$ (see Lemmas \ref{lem-sc} and \ref{RzG}).
The identity of logarithmic decay on \eqref{lwlg} is revealed in our main result.
Asymptotic expansion of this type is firstly introduced by Escobedo and Zuazua \cite{EZ} for the convection-diffusion equation.
On several frameworks, large-time behavior and asymptotic expansion of Navier-Stokes flow are studied by many authors (for example, see \cite{Brndls-Vgnrn,Brndls-Okb,ChJn}).
In \cite{Okb18}, an asymptotic expansion is provided without the moment condition on initial data.
Moreover, by setting the Hardy space as framework, they draw spatial decay of solutions.
Here we choose the other way to lead the spatial structure of the solution, i.e., we study the estimates in weighted spaces as in \cite{Kkvc-Ris}.
To describe far-field asymptotics, we use
\begin{equation}\label{Umt}
	U_m^t (t)
	=
	\sum_{|\beta|=m} \frac{\nabla^\beta R^\bot RG (t)}{\beta!}
	\int_0^t \int_{\mathbb{R}^2} (-y)^\beta (\omega u) (s,y) dyds
\end{equation}
instead of $U_m^\infty$ for $m=1$ and $2$.
\begin{proposition}\label{prop-lows}
	Let $\omega_0 \in L^1 (\mathbb{R}^2) \cap L^\infty (\mathbb{R}^2),~ |x|^5 \omega_0 \in L^1 (\mathbb{R}^2)$ and $\int_{\mathbb{R}^2} x^\alpha \omega_0 (x) dx = 0$ for $|\alpha| \le 1$.
	Assume that the solutions $u$ of \eqref{NS} for $a = -\nabla^\bot (-\Delta)^{-1} \omega_0$ and $\omega$ of \eqref{mcurl} fulfill \eqref{decay} and \eqref{decay-vort} for $k = 5$, respectively.
	Then
	\[
		\Bigl\| |x|^\mu \Bigl( u(t) - \sum_{m=1}^2 (U_m+U_m^t) (t) \Bigr) \Bigr\|_{L^1 (\mathbb{R}^2)}
		= O ( t^{-\frac32 + \frac\mu{2}} \log t )
	\]
	as $t \to +\infty$ for $0 \le \mu < 3$, and
	\[
		\Bigl\| |x|^\mu \Bigl( u(t) - \sum_{m=1}^2 (U_m+U_m^t) (t) \Bigr) \Bigr\|_{L^\infty (\mathbb{R}^2)}
		= O ( t^{-\frac52 + \frac\mu{2}} \log t )
	\]
	as $t \to +\infty$ for $0 \le \mu \le 5$.	
\end{proposition}
This proposition is shown in the similar way as in proof of our main result (see the sentences after proof of Theorem \ref{thm-st}).
From the same view point as above, we see from this proposition that $U_m + U_m^t$ gives far-field asymptotics of $u$.
Indeed
$|x|^5 (U_m + U_m^t) \not\in L^\infty (\mathbb{R}^2)$.
Moreover, large-time behavior of $U_m^t$ also is clear.
In fact, applying \eqref{decay} and \eqref{decay-vort} to the right-hand side of
$
	U_m^t (t) - U_m^\infty (t)
	=
	- \sum_{|\beta|=m} \frac{\nabla^\beta R^\bot RG (t)}{\beta!} \int_t^\infty \int_{\mathbb{R}^2} (-y)^\beta (\omega u) (s,y) dyds,
$
we obtain that $\| U_m^t (t) - U_m^\infty (t) \|_{L^q (\mathbb{R}^2)} \le Ct^{-\gamma_q-\frac{m}2} (1+t)^{-\frac32 + \frac{m}2}$ for $m = 1$ and $2$, and $1 \le q \le \infty$.
Namely $U_m^t$ converges to $U_m^\infty$ as $t \to +\infty$ asymptotically.
Those profiles are obtained by the following procedure.
By expanding the nonlinear term on \eqref{msu}, we see that
\[
	\int_0^t R^\bot RG (t-s) * (\omega u) (s) ds
	=
	\sum_{m=1}^2 U_m^t (t) + r_2 (t),
\]
where
\[
	r_2 (t) = \int_0^t \int_{\mathbb{R}^2} \Bigl(
		R^\bot RG (t-s,x-y) - \sum_{2l+|\beta|=0}^2 \frac{\partial_t^l \nabla^\beta R^\bot R G (t,x)}{\beta!} (-s)^l (-y)^\beta
	\Bigr) (\omega u) (s,y) dyds.
\]
Here we used that $\int_{\mathbb{R}^2} (\omega u) (s,y) dy = 0$.
From the term of initial-data, we get $U_m$ (see Lemmas \ref{lem-asymplint} and \ref{lem-asymplinsp}).
Taylor theorem guarantees that the remained term $r_2$ decays fast and we conclude Propositions \ref{prop-lowt} and \ref{prop-lows}.
Those two propositions give the asymptotic expansions of $u$ with second order.
The renormalization in time yields one with higher-order.
For some related equations, the theory for renormalization is developed in \cite{Iwbch,KtM,Ngi-Ymd}.
To apply the renormalization to our model, asymptotic profiles of $\omega$ are required.
From Biot-Savard law, it is natural that the profiles of $\omega$ are given by $\Omega_m = \nabla \times (U_{m-1} + U_{m-1}^\infty)$, i.e.,
\begin{equation}\label{defOmg}
	\Omega_m (t) = \sum_{|\alpha|=m} \frac{\nabla^\alpha G(t)}{\alpha!} \int_{\mathbb{R}^2} (-y)^\alpha \omega_0 (y) dy
	- \sum_{|\beta|=m-1} \frac{\nabla^\beta \nabla G(t)}{\beta!} \cdot \int_0^\infty \int_{\mathbb{R}^2} (-y)^\beta (\omega u) (s,y) dyds
\end{equation}
for $m = 2$ and $3$.
Those terms are also derived from \eqref{mcurl} directly through the same way as in \cite{Crpo,EZ,Fjgk-Mykw}.
Then the asymptotic profiles of $\omega u$ on \eqref{msu} are provided by the products of $\Omega_m$ and $U_m + U_m^\infty$.
More precisely, $\omega u$ converges to the sum of
\begin{equation}\label{defI}
	\mathcal{I}_5 (t)
	=
	\Omega_2 (U_1 + U_1^\infty) (t)
	\quad
	\text{and}
	\quad
	\mathcal{I}_6 (t)
	=
	\sum_{m=1}^2 \Omega_{4-m} (U_m + U_m^\infty) (t)
\end{equation}
asymptotically (see Corollary \ref{cor-asympf}).
Throughout this paper, the indexes under symbols indicate their scalings.
Namely $\lambda^{2+m} \Omega_m (\lambda^2 t, \lambda x) = \Omega_m (t,x)$ and
\begin{equation}\label{scI}
	\lambda^{2+p} \mathcal{I}_p (\lambda^2 t, \lambda x) = \mathcal{I}_p (t,x)
\end{equation}
hold for $\lambda > 0$.
Clearly, those functions satisfy that $|x|^\mu \Omega_m (1),~ |x|^\mu \mathcal{I}_p (1) \in L^\infty (\mathbb{R}^2)$ for any $\mu \ge 0$.
%
%
By renormalizing $\omega u$ by $\mathcal{I}_5$, we expand the above remained term $r_2$:
\[
\begin{split}
	r_2 (t)
	&=
	\sum_{2l+|\beta|=3} \frac{\partial_t^l\nabla^\beta R^\bot R G(t)}{\beta!}
	\int_0^t \int_{\mathbb{R}^2} (-s)^l (-y)^\beta (\omega u) (s,y) dyds\\
	&+ \int_0^t \int_{\mathbb{R}^2}
		\Bigl(
			R^\bot RG(t-s,x-y) - \sum_{2l+|\beta|=0}^3 \frac{\partial_t^l \nabla^\beta R^\bot RG(t,x)}{\beta!} (-s)^l (-y)^\beta
		\Bigr) (\omega u) (s,y)
	dyds\\
	&=
	U_3^t (t) + K_3 (t) + J_3 (t) + r_3 (t),
\end{split}
\]
where
\begin{equation}\label{U3t}
\begin{split}
	U_3^t (t) &= \sum_{2l+|\beta|=3} \frac{\partial_t^l \nabla^\beta R^\bot R G(t)}{\beta!}
	\int_0^t \int_{\mathbb{R}^2} (-s)^l (-y)^\beta \left( (\omega u) (s,y) - \mathcal{I}_5 (1+s,y) \right) dyds,\\
	K_3 (t) &= \sum_{2l+|\beta|=3} \frac{\partial_t^l \nabla^\beta R^\bot RG(t)}{\beta!}
	\int_0^t \int_{\mathbb{R}^2} (-s)^l (-y)^\beta \mathcal{I}_5 (1+s,y) dyds,\\
	J_3 (t) &= \int_0^t \int_{\mathbb{R}^2}
		\Bigl(
			R^\bot RG(t-s,x-y) - \sum_{2l+|\beta|=0}^3 \frac{\partial_t^l \nabla^\beta R^\bot RG(t,x)}{\beta!} (-s)^l (-y)^\beta
		\Bigr) \mathcal{I}_5 (s,y)
	dyds
\end{split}
\end{equation}
and
\[
\begin{split}
	r_3 (t)
	&= \int_0^t \int_{\mathbb{R}^2}
		\Bigl(
			R^\bot RG(t-s,x-y) - \sum_{2l+|\beta|=0}^3 \frac{\partial_t^l \nabla^\beta R^\bot RG(t,x)}{\beta!} (-s)^l (-y)^\beta
		\Bigr)
		\left( \omega u  - \mathcal{I}_5 \right) (s,y)
	dyds.
\end{split}
\]
Now we remark that, if we put $\mathcal{I}_5 (s)$ in $K_3$ instead of $\mathcal{I}_5 (1+s)$, then this term diverges to infinity.
Indeed, from \eqref{scI},
\[
	\int_0^t \int_{\mathbb{R}^2} (-s)^l (-y)^\beta \mathcal{I}_5 (s,y) dyds
	=
	\int_0^t s^{-1} ds \int_{\mathbb{R}^2} (-1)^l (-y)^\beta \mathcal{I}_5 (1,y) dy
\]
since $2l + |\beta| = 3$.
From the same view point, we should confirm that $J_3$ is well-defined (see Proposition \ref{prop-defJ}).
The last step is on the same way, i.e., we renormalize $\omega u - \mathcal{I}_5$ in $r_3$ by $\mathcal{I}_6$:
\[
\begin{split}
	r_3 (t)
	&=
	\sum_{2l+|\beta|=4} \frac{\partial_t^l\nabla^\beta R^\bot R G(t)}{l!\beta!}
	\int_0^t \int_{\mathbb{R}^2} (-s)^l (-y)^\beta \left( \omega u - \mathcal{I}_5 \right) (s,y) dyds\\
	&+ \int_0^t \int_{\mathbb{R}^2}
	\Bigl(
		R^\bot R G(t-s,x-y) - \sum_{2l+|\beta|=0}^4 \frac{\partial_t^l \nabla^\beta R^\bot R G (t,x)}{l!\beta!} (-s)^l (-y)^\beta
	\Bigr) \left( \omega u - \mathcal{I}_5 \right) (s,y) dyds\\
	&= U_4^t (t) + K_4 (t) + J_4 (t) + r_4 (t),
\end{split}
\]
where
\begin{equation}\label{U4t}
\begin{split}
	U_4^t (t) &=
	\sum_{2l+|\beta|=4} \frac{\partial_t^l \nabla^\beta R^\bot R G(t)}{l!\beta!} \int_0^t \int_{\mathbb{R}^2}
		(-s)^l (-y)^\beta \left( (\omega u - \mathcal{I}_5) (s,y) - \mathcal{I}_6 (1+s,y) \right)
	dyds,\\
	K_4 (t) &=
	\sum_{2l+|\beta|=4} \frac{\partial_t^l \nabla^\beta R^\bot R G(t)}{l!\beta!} \int_0^t \int_{\mathbb{R}^2}
		(-s)^l (-y)^\beta \mathcal{I}_6 (1+s,y)
	dyds,\\
	J_4 (t) &=
	\int_0^t \int_{\mathbb{R}^2}
		\Bigl(
			R^\bot RG(t-s,x-y) - \sum_{2l+|\beta|=0}^4 \frac{\partial_t^l \nabla^\beta R^\bot RG(t,x)}{l!\beta!} (-s)^l (-y)^\beta
		\Bigr) \mathcal{I}_6 (s,y)
	dyds
\end{split}
\end{equation}
and
\[
\begin{split}
	r_4 (t) &=
	\int_0^t \int_{\mathbb{R}^2}
		\Bigl(
			R^\bot RG(t-s,x-y)\\
			&\hspace{15mm} - \sum_{2l+|\beta|=0}^4 \frac{\partial_t^l \nabla^\beta R^\bot RG(t,x)}{l!\beta!} (-s)^l (-y)^\beta
		\Bigr)
		\left( \omega u  - \mathcal{I}_5 - \mathcal{I}_6 \right) (s,y)
	dyds.
\end{split}
\]
Combining the above formulas, the nonlinear term of \eqref{msu} is expanded as
\begin{equation}\label{expnl}
	\int_0^t R^\bot R G (t-s) * (\omega u) (s) ds
	=
	\sum_{m=1}^4 U_m^t (t) + \sum_{m=3}^4 (K_m + J_m)(t) + r_4 (t).
\end{equation}
Here the symbol $r$ means the `{\bf r}emained' term.
In fact, we will show that
\begin{equation}\label{r4st}
	\left\| |x|^\mu r_4 (t) \right\|_1 = O ( t^{-\frac52+\frac\mu{2}} (\log t)^2 )
\end{equation}
as $t \to +\infty$ for $0 \le \mu < 5$, and
\begin{equation}\label{r4stinf}
	\left\| |x|^\mu r_4 (t) \right\|_\infty = O ( t^{-\frac72 + \frac{\mu}2} (\log t)^2 )
\end{equation}
 as $t \to +\infty$ for $0 \le \mu \le 7$.
Namely $r_4$ decays fast in far-field.
Since $\omega u$ also decays fast in far-field, renormalization with space-variable is not required.
Then we know that the above procedure is one with time-variable. 
A combination of \eqref{r4st}, \eqref{r4stinf} and Lemma \ref{lem-asymplinsp} in Section 2 provides our main result.
\begin{theorem}\label{thm-st}
Let $\omega_0 \in L^1 (\mathbb{R}^2) \cap L^\infty (\mathbb{R}^2),~ |x|^7 \omega_0 \in L^1 (\mathbb{R}^2)$ and $\int_{\mathbb{R}^2} x^\alpha \omega_0 (x) dx = 0$ for $|\alpha| \le 1$.
Assume that the solutions $u$ of \eqref{NS} for $a = -\nabla^\bot (-\Delta)^{-1} \omega_0$ and $\omega$ of \eqref{mcurl} fulfill \eqref{decay} and \eqref{decay-vort} for $k = 7$, respectively.
Then
\[
	\Bigl\| |x|^\mu \Bigl( u(t) - \sum_{m=1}^4 (U_m + U_m^t) (t) - \sum_{m=3}^4 (K_m + J_m)(t) \Bigr) \Bigr\|_{L^1 (\mathbb{R}^2)}
	=
	O (t^{-\frac{5}2 + \frac\mu2} (\log t)^2)
\]
 as $t \to +\infty$ for $0 \le \mu < 5$, and
\[
	\Bigl\| |x|^\mu \Bigl( u(t) - \sum_{m=1}^4 (U_m + U_m^t) (t) - \sum_{m=3}^4 (K_m + J_m)(t) \Bigr) \Bigr\|_{L^\infty (\mathbb{R}^2)}
	=
	O (t^{-\frac{7}2 + \frac\mu{2}} (\log t)^2 )
\]
 as $t \to +\infty$ for $0 \le \mu \le 7$, where $U_m, U_m^t, K_m$ and $J_m$ are defined by \eqref{Um}, \eqref{Umt}, \eqref{U3t} and \eqref{U4t}.
\end{theorem}
Here we remark that $U_m$ in \eqref{Um} are defined for any $m \in \mathbb{N}$ if $|x|^{m+1} \omega_0 \in L^1 (\mathbb{R}^2)$.
This theorem suggests that the remained ingredients of velocity are decaying or growing slowly in far-field as $t \to +\infty$.
The functions on expansion have the following structures.
Firstly, H\"armander-Mikhlin type estimate says that $|x|^7 U_m, |x|^7 U_m^t, |x|^7 K_m$ are not in $L^\infty (\mathbb{R}^2)$ (see Lemma \ref{RzG}).
For $\lambda > 0$, $U_m$ and $J_m$ satisfy
\begin{equation}\label{scU}
	\lambda^{2+m} U_m (\lambda^2 t,\lambda x) = U_m (t,x)
\end{equation}
for $1 \le m \le 4$ and
\begin{equation}\label{scJ}
	\lambda^{2+m} J_m (\lambda^2 t,\lambda x) = J_m (t,x)
\end{equation}
for $m = 3$ and $4$.
For $1 \le q \le \infty$, $U_m^t$ and $K_m$ fulfill that
\begin{equation}\label{dcU}
	\| U_m^t (t) \|_{L^q (\mathbb{R}^2)} = O ( t^{-\gamma_q-\frac{m}2} )
\end{equation}
as $t \to +\infty$ for $1 \le m \le 4$, and
\begin{equation}\label{dcK}
	\| K_m (t) \|_{L^q (\mathbb{R}^2)} = O ( t^{-\gamma_q-\frac{m}2} \log t )
\end{equation}
as $t \to +\infty$ for $m = 3$ and $4$.

The assertions of Theorem \ref{thm-st} with $\mu = 0$ provide sure the large-time behavior of solution.
Indeed, from \eqref{scU} and \eqref{scJ}, we see that $\| U_m (t) \|_{L^q (\mathbb{R}^2)} = t^{-\gamma_q-\frac{m}2} \| U_m (1) \|_{L^q (\mathbb{R}^2)}$ and $\| J_m (t) \|_{L^q (\mathbb{R}^2)} = t^{-\gamma_q-\frac{m}2} \| J_m (1) \|_{L^q (\mathbb{R}^2)}$.
Moreover \eqref{dcK} are sharp since $K_m$ are converted to
\[
\begin{split}
	K_3 (t) &= \sum_{|\beta| = 3} \frac{\nabla^\beta R^\bot R G(t)}{\beta!} \int_0^t (1+s)^{-1} ds \int_{\mathbb{R}^2} (-y)^\beta \mathcal{I}_5 (1,y) dy\\
	&+
	\sum_{|\beta| = 1} \partial_t \nabla^\beta R^\bot R G(t) \int_0^t s (1+s)^{-2} ds \int_{\mathbb{R}^2} (-1) (-y)^\beta \mathcal{I}_5 (1,y) dy
\end{split}
\]
and
\[
\begin{split}
	K_4 (t) &= \sum_{|\beta| = 4} \frac{\nabla^\beta R^\bot RG (t)}{\beta!} \int_0^t (1+s)^{-1} ds \int_{\mathbb{R}^2} (-y)^\beta \mathcal{I}_6 (1,y) dy\\
	&+
	\sum_{|\beta| = 2} \frac{\partial_t \nabla^\beta R^\bot R G(t)}{\beta!} \int_0^t s (1+s)^{-2} ds \int_{\mathbb{R}^2} (-1) (-y)^\beta \mathcal{I}_6 (1,y) dy.
\end{split}
\]
Here we used \eqref{scI} and that $\int_{\mathbb{R}^2} \mathcal{I}_6 (1,y) dy = 0$.
However, as far as we concern \eqref{dcU}, behaviors of $U_m^t$ as $t \to +\infty$ are not clear.
In order to clarify it, we introduce the better functions.
For $1 \le m \le 3$, we expand $U_m^t$, then
\[
	U_1^t (t) + U_2^t (t) + U_3^t (t)
	=
	U_1^\infty (t) + U_2^\infty (t) + U_3^\infty (t) + V_3 (t) + V_4 (t) + r_4^\infty (t),
\]
where $U_1^\infty$ and $U_2^\infty$ are defined by \eqref{Um}, and
\begin{equation}\label{U3inf}
\begin{split}
	U_3^\infty (t)
	&= \sum_{2l+|\beta|=3} \frac{\partial_t^l \nabla^\beta R^\bot R G(t)}{\beta!}
	\int_0^\infty \int_{\mathbb{R}^2} (-s)^l (-y)^\beta \left( \omega u (s,y) - \mathcal{I}_5 (1+s,y) \right) dyds,\\
	V_3 (t)
	&= -\sum_{|\beta| = 1}^2 \frac{\nabla^\beta R^\bot RG(t)}{\beta!}
	\int_t^\infty \int_{\mathbb{R}^2} (-y)^\beta \mathcal{I}_5 (s,y) dyds,\\
	V_4 (t)
	&= -\sum_{2l + |\beta| = 1}^3 \frac{\partial_t^l \nabla^\beta R^\bot RG(t)}{\beta!}
	\int_t^\infty \int_{\mathbb{R}^2} (-s)^l (-y)^\beta \mathcal{I}_6 (s,y) dyds
\end{split}
\end{equation}
and
\[
\begin{split}
	r_4^\infty (t)
	&=
	\sum_{2l + |\beta|=1}^3 \frac{\partial_t^l \nabla^\beta R^\bot R G(t)}{\beta!} \int_t^\infty \int_{\mathbb{R}^2}
		(-s)^l (-y)^\beta \left( \omega u - \mathcal{I}_5 - \mathcal{I}_6 \right)(s,y)
	dyds\\
	&+
	\sum_{2l+|\beta|=3} \frac{\partial_t^l \nabla^\beta R^\bot R G(t)}{\beta!} \int_t^\infty \int_{\mathbb{R}^2}
		(-s)^l (-y)^\beta \left( \mathcal{I}_5 (s,y) - \mathcal{I}_5 (1+s,y) \right)
	dyds.
\end{split}
\]	
For $m = 4$, we choose
\begin{equation}\label{U4inf}
	U_4^\infty (t) =
	\sum_{2l+|\beta|=4} \frac{\partial_t^l \nabla^\beta R^\bot R G(t)}{l!\beta!} \int_0^\infty \int_{\mathbb{R}^2}
		(-s)^l (-y)^\beta \left( \omega u (s,y) - \mathcal{I}_5 (s,y) - \mathcal{I}_6 (1+s,y) \right)
	dyds
\end{equation}
instead of $U_4^t$.
We confirm later that
\begin{equation}\label{Usa}
	\| U_3^\infty (t) - U_3^t (t) \|_{L^q (\mathbb{R}^2)}
	=
	O (t^{-\gamma_q-2})\quad
	\text{and}\quad
	\| U_4^\infty (t) - U_4^t (t) \|_{L^q (\mathbb{R}^2)}
	=
	O (t^{-\gamma_q-\frac52} \log t)
\end{equation}
as $t \to +\infty$.
%
Here, for $\lambda > 0$, we see that
\begin{equation}\label{scUinf}
	\lambda^{2+m} U_m^\infty (\lambda^2 t, \lambda x) = U_m^\infty (t,x)
\end{equation}
for $1 \le m \le 4$, and that
\begin{equation}\label{scV}
	\lambda^{2+m} V_m (\lambda^2 t, \lambda x) = V_m (t,x)
\end{equation}
for $m = 3$ and $4$.
Moreover we will show that
\begin{equation}\label{r4inft}
	\| r_4^\infty (t) \|_{L^q (\mathbb{R}^2)}
	=
	O (t^{-\gamma_q - \frac52} \log t)
\end{equation}
as $t \to +\infty$ for $1 \le q \le \infty$.
Therefore we conclude \eqref{dcU} and obtain our second assertion.
\begin{theorem}\label{thm-t}
Let $\omega_0 \in L^1 (\mathbb{R}^2) \cap L^\infty (\mathbb{R}^2),~ |x|^5 \omega_0 \in L^1 (\mathbb{R}^2)$ and $\int_{\mathbb{R}^2} x^\alpha \omega_0 (x) dx = 0$ for $|\alpha| \le 1$.
Assume that the solutions $u$ of \eqref{NS} for $a = -\nabla^\bot (-\Delta)^{-1} \omega_0$ and $\omega$ of \eqref{mcurl} fulfill \eqref{decay} and \eqref{decay-vort} for $k = 5$, respectively.
Then
\[
	\Bigl\| u(t) - \sum_{m=1}^4 (U_m + U_m^\infty) (t) - \sum_{m=3}^4 (K_m + J_m + V_m)(t) \Bigr\|_{L^q (\mathbb{R}^2)}
	=
	o (t^{-\gamma_q-2})
\]
as $t \to +\infty$ for $1 \le q \le \infty$, where $U_m, U_m^\infty, K_m, J_m$ and $V_m$ are defined by \eqref{Um}, \eqref{Uminf}, \eqref{U3t}, \eqref{U4t}, \eqref{U3inf} and \eqref{U4inf}
In addition, if $|x|^6 \omega_0 \in L^1 (\mathbb{R}^2)$, then
\[
	\Bigl\| u(t) - \sum_{m=1}^4 (U_m + U_m^\infty) (t) - \sum_{m=3}^4 (K_m + J_m + V_m)(t) \Bigr\|_{L^q (\mathbb{R}^2)}
	=
	O (t^{-\gamma_q-\frac52} (\log t)^2)
\]
as $t \to +\infty$ for $1 \le q \le \infty$.
\end{theorem}
Since \eqref{scU}, \eqref{scJ}, \eqref{dcK}, \eqref{scUinf} and \eqref{scV} hold, large-time behaviors of any terms on the expansion are clear.
We emphasize that the first assertion is sharp under the condition $|x|^5 \omega_0 \in L^1 (\mathbb{R}^2)$ since this assumption is compatible with $|x|^4 a \in L^1 (\mathbb{R}^2)$ (see Lemmas \ref{lem-asymplint} and \ref{lem-asymplinsp} in Section 2).

Theorems \ref{thm-st} and \ref{thm-t} provide the sharp estimates for asymptotic expansions with fourth order.
If we try to develop such estimates based on \eqref{HFK}, then at least the property $|x|^4 u \otimes u \in L_{loc}^1 (0,\infty,L^1 (\mathbb{R}^2))$ is required.
However this contradicts \eqref{sp-decay}.
Hence we employ the vorticity and adopt \eqref{msu}.

\vspace{3mm}

\paragraph{\it Notations.}
For vectors, we abbreviate them by using the same letters, for example, $a = (a_1, a_2)$.
For $x = (x_1,x_2)$ and $y = (y_1,y_2) \in \mathbb{R}^2$, we denote $x \cdot y = x_1 y_1 + x_2 y_2$ and $|x|^2 = x \cdot x$.
We often omit the spatial variable from functions, so $U(t) = U(t,x)$.
For vector-fields $f$ and $g$, the convolution of them is simply denoted by $f*g (x) = \int_{\mathbb{R}^2} f(x-y) \cdot g (y) dy = \int_{\mathbb{R}^2} f(y) \cdot g(x-y) dy$.
Hence $f*g$ is scalar here.
We symbolize that $\partial_t = \partial/\partial t,~ \partial_1 = \partial/\partial x_1,~ \partial_2 = \partial/\partial x_2,~ \nabla = (\partial_1,\partial_2),~ \nabla^\bot = (-\partial_2, \partial_1)$ and $\Delta = |\nabla|^2 = \partial_1^2 + \partial_2^2$.
The length of a multi-index $\alpha = (\alpha_1,\alpha_2) \in \mathbb{Z}_+^2$ is given by $|\alpha| = \alpha_1 + \alpha_2$, where $\mathbb{Z}_+ = \mathbb{N} \cup \{ 0 \}$.
We abbreviate that $\alpha ! = \alpha_1 ! \alpha_2 !,~ x^\alpha = x_1^{\alpha_1} x_2^{\alpha_2}$ and $\nabla^\alpha = \partial_1^{\alpha_1} \partial_2^{\alpha_2}$.
We define the Fourier transform and its inverse by $\hat{\varphi} (\xi) = \mathcal{F} [\varphi] (\xi) = (2\pi)^{-1} \int_{\mathbb{R}^2} \varphi (x) e^{-ix\cdot\xi} dx$ and $\check{\varphi} (x) = \mathcal{F}^{-1} [\varphi] (x) = (2\pi)^{-1} \int_{\mathbb{R}^2} \varphi (\xi) e^{ix\cdot\xi} d\xi$, respectively, where $i = \sqrt{-1}$.
The Riesz transforms are defined by $R_j \varphi = \partial_j (-\Delta)^{-1/2} \varphi = \mathcal{F}^{-1} [\frac{i\xi_j}{|\xi|} \hat{\varphi}]$ for $j = 1$ and $2$.
For the tensor $R^\bot R$, see Definition \ref{def-Rz2} and the added sentence.
Similarly $\nabla^\bot (-\Delta)^{-1} \varphi = \mathcal{F}^{-1} [i\xi^\bot |\xi|^{-2} \hat{\varphi}]$, where $\xi^\bot = (-\xi_2, \xi_1)$.
For $1 \le q \le \infty$, $L^q (\mathbb{R}^n)$ denotes the Lebesgue space and $\| \cdot \|_{L^q (\mathbb{R}^n)}$ is its norm.
For a vector-field $f = f(x)$ and a tensor $F = F(x)$, we abbreviate their norms as $\| f \|_{L^q (\mathbb{R}^2)} = \| f \|_{L^q (\mathbb{R}^2, \mathbb{R}^2)}$ and $\| F \|_{L^q (\mathbb{R}^2)} = \| F \|_{L^q (\mathbb{R}^2, \mathbb{R}^{2\times 2})}$.
We denote the two dimensional gaussian and its decay-rate in $L^q (\mathbb{R}^2)$ by $G(t,x) = (4\pi t)^{-1} e^{-|x|^2/4t}$ and $\gamma_q = 1 - \frac1q$, respectly.
Namely $\| G(t) \|_{L^q (\mathbb{R}^2)} = t^{-\gamma_q} \| G(1) \|_{L^q (\mathbb{R}^2)}$ for $t > 0$.
Throughout this paper, the indexes under symbols mean those scalings or decay-rates in time.
For example, $\lambda^{2+m} J_m (\lambda^2 t, \lambda x) = J_m (t,x)$ for $\lambda > 0$, and $\| K_m (t) \|_{L^q (\mathbb{R}^2)} = O (t^{-\gamma_q - \frac{m}2} \log t)$ as $t \to +\infty$ for $1 \le q \le \infty$.
Various positive constants are simply denoted by $C$.

\section{Preliminaries}
In this section, we prepare several lemmas which are used to confirm our main results.
Structures of the asymptotic expansions are clarified by using the following lemmas.
\begin{lemma}\label{lem-sc}
Let $1 \le q \le \infty,~ m, \mu \ge 0$ and a measurable function $U = U(t,x)$ fulfill that $|x|^\mu U(1) \in L^q (\mathbb{R}^2)$ and $\lambda^{2+m} U(\lambda^2 t, \lambda x) = U(t,x)$ for $\lambda > 0$.
Then $\| |x|^\mu U(t) \|_{L^q (\mathbb{R}^2)} = t^{-\gamma_q-\frac{m}2+\frac\mu{2}} \| |x|^\mu U(1) \|_{L^q (\mathbb{R}^2)}$.
\end{lemma}
%
%
\begin{lemma}\label{RzG}
For $\alpha \in \mathbb{Z}_+^2,~ (1+|x|)^{1+|\alpha|} \nabla^\alpha \nabla^\bot (-\Delta)^{-1} G (1)$ and $(1+|x|)^{2+|\alpha|} \nabla^\alpha R^\bot R G(1)$ are bounded, and $\lambda^{1+|\alpha|} \nabla^\alpha \nabla^\bot (-\Delta)^{-1} G (\lambda^2 t, \lambda x) = \nabla^\alpha \nabla^\bot (-\Delta)^{-1} G (t,x)$ and $\lambda^{2+|\alpha|} \nabla^\alpha R^\bot R G(\lambda^2 t, \lambda x) = \nabla^\alpha R^\bot R$ $G(t,x)$ hold for $\lambda > 0$.
\end{lemma}
Proof of Lemma \ref{lem-sc} is straightforward and Lemma \ref{RzG} is confirmed by H\"ormander-Mikhlin estimate (see for example \cite{Stin,Zmr} and \cite[Theorem 2.3]{Sbt-Smz}) and elemantary calculus.

We consider the linear heat equation with initial-data $a = -\nabla^\bot (-\Delta)^{-1} \omega_0$, then we see the following lemmas.
\begin{lemma}\label{lem-asymplint}
Let  $\omega_0 \in L^1 (\mathbb{R}^2) \cap L^\infty (\mathbb{R}^2),~ |x|^5 \omega_0 \in L^1 (\mathbb{R}^2)$ and $\int_{\mathbb{R}^2} x^\alpha \omega_0 (x) dx = 0$ for $|\alpha| \le 1$.
Then
\[
		\Bigl\| \nabla^\bot (-\Delta)^{-1} G(t) * \omega_0 + \sum_{m=1}^4 U_m (t) \Bigr\|_{L^q (\mathbb{R}^2)}
		= o ( t^{-\gamma_q-2} )
\]
as $t \to +\infty$ for $1 \le q \le \infty$ , where $U_m$ are defined by \eqref{Um}.
In addition, if $|x|^6 \omega_0 \in L^1 (\mathbb{R}^2)$, then
\[
		\Bigl\| \nabla^\bot (-\Delta)^{-1} G(t) * \omega_0 + \sum_{m=1}^4 U_m (t) \Bigr\|_{L^q (\mathbb{R}^2)}
		= O ( t^{-\gamma_q-\frac52} )
\]
as $t \to +\infty$ for $1 \le q \le \infty$.
\end{lemma}
\begin{proof}
We separate the domain to $\mathbb{R}^2 = \{ |y| \le R(t) \} \cup \{ |y| > R(t) \}$ by a positive function $R(t)$ such that $R(t) = o (t^{1/2})$ as $t \to +\infty$, then by Taylor theorem we see that
\[
\begin{split}
	&\nabla^\bot (-\Delta)^{-1} G(t) * \omega_0 + \sum_{m=1}^4 U_m (t)
	=
	\sum_{|\alpha|=6} \int_{|y| \le R(t)} \int_0^1
		\frac{\nabla^\alpha \nabla^\bot (-\Delta)^{-1} G(t,x-\lambda y)}{\alpha!} \lambda^5 (-y)^\alpha \omega_0 (y)
	d\lambda dy\\
	&+
	\sum_{|\alpha|=5} \int_{|y| > R(t)} \int_0^1
		\frac{\nabla^\alpha \nabla^\bot (-\Delta)^{-1} G(t,x-\lambda y)}{\alpha!} \lambda^4 (-y)^\alpha \omega_0 (y)
	d\lambda dy\\
	&-
	\sum_{|\alpha|=5} \frac{\nabla^\alpha \nabla^\bot (-\Delta)^{-1} G(t)}{\alpha!} \int_{|y| > R(t)} 
		(-y)^\alpha \omega_0 (y)
	dy.
\end{split}
\]
%
Thus, by Lemmas \ref{lem-sc} and \ref{RzG}, we get that
\[
\begin{split}
	&\Bigl\| \nabla^\bot (-\Delta)^{-1} G(t) * \omega_0 + \sum_{m=1}^4 U_m (t) \Bigr\|_{L^q (\mathbb{R}^2)}\\
	&\le
	C t^{-\gamma_q - \frac52} R(t) \| |x|^5 \omega_0 \|_{L^1 (\mathbb{R}^2)}
	+
	C t^{-\gamma_q-2} \int_{|y| > R(t)} |y|^5 |\omega_0 (y)| dy
\end{split}
\]
and conclude the first assertion.
Proof for the second assertion is straightforward.
\end{proof}
Here we used `Hausdorff-Young-type inequality' in the last estimate.
Namely, for $f \in L^1 (\mathbb{R}^2)$ and $g \in L^q (\mathbb{R}^2)$, we see from H\"older inequality that
$| \int_{\mathbb{R}^2} f(y) \int_0^1 g(x-\lambda y) d\lambda dy |
\le \| f \|_{L^1 (\mathbb{R}^2)}^{1-1/q} (\int_{\mathbb{R}^2} |f(y)| \int_0^1 |g(x-\lambda y)|^q d\lambda dy)^{1/q}$.
Thus Fubini theorem says that $\| \int_{\mathbb{R}^2} f(y) \int_0^1 g(x-\lambda y) d\lambda dy \|_{L^q (\mathbb{R}^2)} \le \| f \|_{L^1 (\mathbb{R}^2)} \| g \|_{L^q (\mathbb{R}^2)}$.
Hereafter we call the estimates of this type also by `Hausdorff-Young inequality' simply.

\begin{lemma}\label{lem-asymplinsp}
Let  $\omega_0 \in L^1 (\mathbb{R}^2) \cap L^\infty (\mathbb{R}^2),~ |x|^7 \omega_0 \in L^1 (\mathbb{R}^2)$ and $\int_{\mathbb{R}^2} x^\alpha \omega_0 (x) dx = 0$ for $|\alpha| \le 1$.
Then
	\[
		\Bigl\| |x|^\mu \Bigl( \nabla^\bot (-\Delta)^{-1} G(t) * \omega_0 + \sum_{m=1}^4 U_m (t) \Bigr) \Bigr\|_{L^1 (\mathbb{R}^2)}
		= O ( t^{-\frac52 + \frac\mu2} )
	\]
	as $t \to +\infty$ for $0 \le \mu < 5$, and
	\[
		\Bigl\| |x|^\mu \Bigl( \nabla^\bot (-\Delta)^{-1} G(t) * \omega_0 + \sum_{m=1}^4 U_m (t) \Bigr) \Bigr\|_{L^\infty (\mathbb{R}^2)}
		= O ( t^{-\frac72 + \frac{\mu}2} )
	\]
as $t \to +\infty$ for $0 \le \mu \le 7$.
\end{lemma}
\begin{proof}
Taylor theorem also gives that
\[
	\nabla^\bot (-\Delta)^{-1} G(t) * \omega_0 + \sum_{m=1}^4 U_m (t)
	=
	\rho_4^1 (t) + \rho_4^2 (t),
\]
where
\[
\begin{split}
	\rho_4^1 (t)
	&=
	\sum_{|\alpha|=6} \int_{|y| \le |x|/2} \int_0^1
		\frac{\nabla^\alpha \nabla^\bot (-\Delta)^{-1} G(t,x-\lambda y)}{\alpha!} \lambda^5 (-y)^\alpha \omega_0 (y)
	d\lambda dy
\end{split}
\]
and
\[
\begin{split}
	\rho_4^2 (t)
	&=
	\sum_{|\alpha|=2} \int_{|y| > |x|/2} \int_0^1
		\frac{\nabla^\alpha \nabla^\bot (-\Delta)^{-1} G(t,x-\lambda y)}{\alpha!} \lambda (-y)^\alpha \omega_0 (y)
	d\lambda dy\\
	&-
	\sum_{|\alpha|=2}^5 \frac{\nabla^\alpha \nabla^\bot (-\Delta)^{-1} G(t)}{\alpha!} \int_{|y| > |x|/2} 
		(-y)^\alpha \omega_0 (y)
	dy.
\end{split}
\]
Hence
\[
\begin{split}
	\| |x|^\mu \rho_4^1 (t) \|_{L^1 (\mathbb{R}^2)}
	&\le
	C \sum_{|\alpha| = 6} \| |x|^\mu \nabla^\alpha \nabla^\bot (-\Delta)^{-1} G(t) \|_{L^1 (\mathbb{R}^2)}
	\left\| x^\alpha \omega \right\|_{L^1 (\mathbb{R}^2)}
\end{split}
\]
and
\[
\begin{split}
	\| |x|^\mu \rho_4^2 (t) \|_{L^1 (\mathbb{R}^2)}
	&\le
	C \sum_{|\alpha|=2}^5 \| |x|^{|\alpha|-2} \nabla^\alpha \nabla^\bot (-\Delta)^{-1} G(t) \|_{L^1 (\mathbb{R}^2)}
	\| |x|^{\mu - |\alpha| + 2} x^\alpha \omega_0 \|_{L^1 (\mathbb{R}^2)}.
\end{split}
\]
Thus we get the first assertion from Lemmas \ref{lem-sc} and \ref{RzG}.
Here the second part originally is
\[
\begin{split}
	\rho_4^2 (t)
	&=
	\int_{|y| > |x|/2}
		\nabla^\bot (-\Delta)^{-1} G(t,x- y) \omega_0 (y)
	dy\\
	&-
	\sum_{|\alpha|=2}^5 \frac{\nabla^\alpha \nabla^\bot (-\Delta)^{-1} G(t)}{\alpha!} \int_{|y| > |x|/2} 
		(-y)^\alpha \omega_0 (y)
	dy.
\end{split}
\]
Therefore, in the same way as above, we obtain the second assertion, and complete the proof.
\end{proof}
By the way, upon the condition on above lemmas, we see that $\| \nabla^\bot (-\Delta)^{-1} G(t) * \omega_0 \|_{L^1 (\mathbb{R}^2)} \le C t^{-1/2}$ since Taylor theorem gives that
\[
	\nabla^\bot (-\Delta)^{-1} G(t) * \omega_0
	=
	\sum_{|\alpha| = 2} \int_{\mathbb{R}^2} \int_0^1
		\frac{\nabla^\alpha \nabla^\bot (-\Delta)^{-1} G(t,x-\lambda y)}{\alpha!}
		\lambda (-y)^\alpha \omega_0 (y)
	d\lambda dy.
\]
Therefore, in \eqref{decay}, case $q = 1$ requires special treatment.
However it is not essential in our proofs.

In order to employ the renormalization, we prepare the following estimate:
%
%
%
%
%
%
%
\begin{proposition}\label{prop-AsympVort}
Let $\omega_0 \in L^1 (\mathbb{R}^2) \cap L^\infty (\mathbb{R}^2),~ |x|^4 \omega_0 \in L^1 (\mathbb{R}^2)$ and $\int_{\mathbb{R}^2} x^\alpha \omega_0 (x) dx = 0$ for $|\alpha| \le 1$.
Assume that the solutions $u$ of \eqref{NS} for $a = -\nabla^\bot (-\Delta)^{-1} \omega_0$ and $\omega$ of \eqref{mcurl} fulfill \eqref{decay} and \eqref{decay-vort} for $k = 3$.
Then $\omega = \nabla \times u$ satisfies that
\begin{equation}\label{AsympVort}
	\| \omega (t) - \Omega_2 (t) - \Omega_3 (t) \|_{L^q (\mathbb{R}^2)}
	\le
	C t^{-\gamma_q-\frac32} (1+t)^{-\frac12} \log (2+t)
\end{equation}
for $1 \le q \le \infty$, where $\Omega_m$ are defined by \eqref{defOmg}.
In addition, if $|x|^k \omega_0 \in L^1 (\mathbb{R}^2)$ and $\omega$ satisfies \eqref{decay-vort} for some $k \ge 4$, then
\begin{equation}\label{AsympVortWt}
	\left\| |x|^k ( \omega - \Omega_2 (t) - \Omega_3 (t) ) \right\|_{L^q (\mathbb{R}^2)}
	\le
	C t^{-\gamma_q} (1+t)^{-2+\frac{k}{2}} \log (2+t)
\end{equation}
holds for $1 \le q \le \infty$.
\end{proposition}
\begin{proof}
The first statement \eqref{AsympVort} is derived from the same procedure as in \cite{Crpo,EZ,Fjgk-Mykw} together with \eqref{decay-vort}.
The derivation process of \eqref{AsympVortWt} as following without $|x|^k$ also gives \eqref{AsympVort}.
Hence we show only \eqref{AsympVortWt}.
Firstly \eqref{decay-vort} and the scaling of $\Omega_m$ immediately give that
$
	\| |x|^k ( \omega - \Omega_2 (t) - \Omega_3 (t) )\|_{L^q (\mathbb{R}^2)}
	\le
	C t^{-\gamma_q} (1+t)^{-1+\frac{k}{2}}.
$
Hence the singularity as $t \to +0$ is bouded by $C t^{-\gamma_q}$.
Since $\int_{\mathbb{R}^2} x^\alpha \omega_0 (x) dx = 0$ for $|\alpha| \le 1$ and $\int_{\mathbb{R}^2} (\omega u) (t,x) dx = 0$, we see that
\[
\begin{split}
	&\omega (t) - \Omega_2 (t) - \Omega_3 (t)
	= \int_{\mathbb{R}^2} \biggl( G(t,x-y) - \sum_{|\alpha| = 0}^3 \frac{\nabla^\alpha G(t)}{\alpha!} (-y)^\alpha \biggr) \omega_0 (y) dy\\
	&-
	\int_0^t \int_{\mathbb{R}^2} \biggl(
		\nabla G (t-s,x-y) - \sum_{2l+|\beta|=0}^2 \frac{\partial_t^l \nabla^\beta \nabla G(t,x)}{\beta!} (-s)^l (-y)^\beta
	\biggr) \cdot (\omega u) (s,y)
	dyds\\
	&+
	\sum_{2l+|\beta| = 1}^2 \frac{\partial_t^l \nabla^\beta \nabla G(t)}{\beta!}
	\cdot \int_t^\infty \int_{\mathbb{R}^2} (-s)^l (-y)^\beta (\omega u) (s,y) dyds.
\end{split}
\]
The estimate for the first term on right-hand side is well-known.
For the last term, Lemmas \ref{lem-sc} and \ref{RzG}, \eqref{decay} and \eqref{decay-vort} yield that
\[
	\biggl\| |x|^k \sum_{2l+|\beta| = 1}^2 \frac{\partial_t^l \nabla^\beta \nabla G(t)}{\beta!} \cdot \int_t^\infty \int_{\mathbb{R}^2} (-s)^l (-y)^\beta (\omega u) (s,y) dyds \biggr\|_{L^q (\mathbb{R}^2)}
	\le
	Ct^{-\gamma_q-\frac32+\frac{k}{2}} (1+t)^{-\frac12}.
\]
In order to estimate the second term, we split the domain $(0,t) \times \mathbb{R}^2$ to $Q_1 \cup Q_2 \cup Q_3$ or $Q_4 \cup Q_5$, where
\begin{gather}
	Q_1 = (0,t/2] \times \{ y \in \mathbb{R}^2~ |~ |y| > |x|/2 \},\quad
	Q_2 = (0,t) \times \{ y \in \mathbb{R}^2~ |~ |y| \le |x|/2 \},\label{Q}\\
	Q_3 = (t/2,t) \times \{ y \in \mathbb{R}^2~ |~ |y| > |x|/2 \},\quad
	Q_4 = Q_2,\quad
	Q_5 = Q_1 \cup Q_3.\notag
\end{gather}
Then
\[
	\int_0^t \int_{\mathbb{R}^2} \biggl(
		\nabla G (t-s,x-y) - \sum_{2l+|\beta|=0}^2 \frac{\partial_t^l \nabla^\beta \nabla G(t,x)}{\beta!} (-s)^l (-y)^\beta
	\biggr) \cdot (\omega u) (s,y)
	dyds
	=
	\varrho_4^1 (t) + \cdots + \varrho_4^5 (t),
\]
where
\[
\begin{split}
	&\varrho_4^j (t)\\
	&=
	\left\{
	\begin{array}{lr}
	\displaystyle
		\iint_{Q_j}
			\biggl( \nabla G(t-s,x-y) - \sum_{l=0}^1 \partial_t^l \nabla G(t,x-y) (-s)^l \biggr)
			\cdot
			(\omega u) (s,y)
		dyds,
		&
		j = 1,2,3,\\
	\displaystyle
		\sum_{l=0}^1 \iint_{Q_j}
			\biggl( \partial_t^l \nabla G(t,x-y) - \sum_{|\beta|=0}^{2-2l} \frac{\partial_t^l \nabla^\beta \nabla G(t,x)}{\beta!} (-y)^\beta \biggr) \cdot
			(-s)^l (\omega u) (s,y)
		dyds,
		&
		j = 4,5.
	\end{array}
	\right.
\end{split}
\]
Taylor theorem leads that
\[
	\varrho_4^1 (t)
	=
	\int_0^{t/2} \int_{|y| > |x|/2} \int_0^1 \frac{\partial_t^2 \nabla G (t-\lambda s,x-y)}{2!} \lambda \cdot (-s)^2 (\omega u) (s,y) d\lambda dyds
\]
and
\[
	\varrho_4^2 (t)
	=
	\int_0^t \int_{|y| \le |x|/2} \int_0^1 \frac{\partial_t^2 \nabla G (t-\lambda s,x-y)}{2!} \lambda \cdot (-s)^2 (\omega u) (s,y) d\lambda dyds.
\]
Thus, from Hausdorff-Young inequality with Lemmas \ref{lem-sc} and \ref{RzG}, \eqref{decay} and \eqref{decay-vort}, we obtain that
\[
	\| |x|^k \varrho_4^1 (t) \|_{L^q (\mathbb{R}^2)}
	\le
	C\int_0^{t/2} \int_0^1
		\| \partial_t^2 \nabla G (t-\lambda s) \|_{L^q (\mathbb{R}^2)} s^2 \| |x|^k (\omega u)(s) \|_{L^1 (\mathbb{R}^2)}
	d\lambda ds
	\le
	C t^{-\gamma_q-\frac52} (1+t)^{\frac12+\frac{k}{2}}
\]
and
\[
\begin{split}
	&\| |x|^k \varrho_4^2 (t) \|_{L^q (\mathbb{R}^2)}
	\le
	C\int_0^{t/2} \int_0^1
		\| |x|^k \partial_t^2 \nabla G (t-\lambda s) \|_{L^q (\mathbb{R}^2)} s^2 \| (\omega u)(s) \|_{L^1 (\mathbb{R}^2)}
	d\lambda ds\\
	&+
	C\int_{t/2}^t \int_0^1
		\| |x|^k \partial_t^2 \nabla G (t-\lambda s) \|_{L^1 (\mathbb{R}^2)} s^2 \| (\omega u)(s) \|_{L^q (\mathbb{R}^2)}
	d\lambda ds
	\le
	C t^{-\gamma_q-\frac52+\frac{k}{2}} (1+t)^{\frac12}.
\end{split}
\]
The estimate for the third term is simpler than the above since Taylor theorem is not required.
In fact
\[
	\| |x|^k \varrho_4^3 (t) \|_{L^q (\mathbb{R}^2)}
	\le
	C t^{-\gamma_q-\frac12} (1+t)^{-\frac32 + \frac{k}{2}}.
\]
By using `Hausdorff-Young inequality' as in the sentence after proof of Lemma \ref{lem-asymplinsp} on the same way as above, we obtain that
\[
	\| |x|^k \varrho_4^4 (t) \|_{L^q (\mathbb{R}^2)}
	\le
	C t^{-\gamma_q-2+\frac{k}{2}} \log (1+t).
\]
The treatment for the last term is straightforward.
Indeed
\[
\begin{split}
	&\| |x|^k \varrho_4^5 (t) \|_{L^q (\mathbb{R}^2)}
	\le
	C\sum_{2l+|\beta|=0}^2 \| x^\beta \partial_t^l \nabla^\beta \nabla G(t) \|_{L^q (\mathbb{R}^2)}
	\int_0^t  s^l \| |x|^k (\omega u) (s) \|_{L^1 (\mathbb{R}^2)} ds\\
	&\le
	C t^{-\gamma_q - \frac12} (1+t)^{-\frac32+\frac{k}{2}}.
\end{split}
\]
Hence we get \eqref{AsympVortWt}.
\end{proof}
\begin{cor}\label{cor-asympf}
	Upon the assumption of Proposition \ref{prop-AsympVort},
	\[
		\left\| |x|^k \left( \omega u - \mathcal{I}_5 - \mathcal{I}_6 \right) (t) \right\|_{L^q (\mathbb{R}^2)}
		\le
		C t^{-\gamma_q} \left( t^{-3+\frac{k}2} + (1+t)^{-3+\frac{k}2} \right) (1+t)^{-\frac12} \log (2+t)
	\]
	holds for $1 \le q \le \infty$, where $\mathcal{I}_5$ and $\mathcal{I}_6$ are defined by \eqref{defI}.
\end{cor}
\begin{proof}
The definitions and elementary calculus provide that
\[
	\omega u - \mathcal{I}_5 - \mathcal{I}_6
	=
	\omega \left( u - U_1 - U_1^\infty - U_2 - U_2^\infty \right)
	+
	\left( \omega - \Omega_2 \right) \left( U_2 + U_2^\infty \right)
	+
	\left( \omega - \Omega_2 - \Omega_3 \right) \left( U_1 + U_1^\infty \right).
\]
By employing \eqref{decay-vort} and Propositions \ref{prop-lowt} and \ref{prop-AsympVort} and the scalings of $U_m$ and $U_m^\infty$ (see also Lemma \ref{lem-sc} and Proposition \ref{prop-defJ}) to the right-hand side, we obtain the decay as $t \to +\infty$.
The singularities as $t \to +0$ are coming from $\omega$ and $\mathcal{I}_6$.
\end{proof}
Of course Proposition \ref{prop-AsympVort} and Corollary \ref{cor-asympf} never give far-field asymptotics of $\omega$ and $\omega u$, respectively, since $|x|^k \Omega_m$ and $|x|^k \mathcal{I}_p$ are in $L^q (\mathbb{R}^2)$.
Before closing this section, we confirm the properties of terms on asymptotic expansions.
\begin{proposition}\label{prop-defJ}
	Under the assumption of Theorem \ref{thm-t}, $U_m, U_m^t, U_m^\infty, K_m, J_m$ and $V_m$ introduced in \eqref{Um}, \eqref{Uminf}, \eqref{Umt}, \eqref{U3t}, \eqref{U4t}, \eqref{U3inf} and \eqref{U4inf} are well-defined on $C((0,\infty), L^1 (\mathbb{R}^2) \cap L^\infty (\mathbb{R}^2))$.
	Moreover, they satisfy \eqref{scU}, \eqref{scJ}, \eqref{dcU}, \eqref{dcK}, \eqref{scUinf} and \eqref{scV}.
\end{proposition}
\begin{proof}
From \eqref{decay}, \eqref{decay-vort} and Proposition \ref{prop-AsympVort}, it is clear that $U_m, U_m^t, U_m^\infty, K_m$ and $V_m$ are in $C((0,\infty),$ $L^1 (\mathbb{R}^2) \cap L^\infty (\mathbb{R}^2))$.
We should show that $J_m$ also are well-defined.
This is done in the same way as in the proof of Proposition \ref{prop-AsympVort}.
We define $J_3^\varepsilon$ for small $\varepsilon > 0$ by
\[
	J_3^\varepsilon (t) =
	\int_\varepsilon^t \int_{\mathbb{R}^2}
		\Bigl(
			R^\bot RG(t-s,x-y) - \sum_{2l+|\beta|=0}^3 \frac{\partial_t^l \nabla^\beta R^\bot RG(t,x)}{\beta!} (-s)^l (-y)^\beta
		\Bigr) \mathcal{I}_5 (s,y)
	dyds,
\]
then, from Taylor theorem, we see that
\[
\begin{split}
	J_3^\varepsilon (t)
	&=
	\int_\varepsilon^{t/2} \int_{\mathbb{R}^2} \int_0^1 \frac{\partial_t^2 R^\bot R G (t-\lambda s,x-y)}{2!} \lambda (-s)^2 \mathcal{I}_5 (s,y) d\lambda dyds\\
	&+
	\sum_{l=0}^1 \sum_{|\beta| = 4-2l} \int_\varepsilon^{t/2} \int_{\mathbb{R}^2} \int_0^1 \frac{\partial_t^l \nabla^\beta R^\bot RG(t,x-\lambda y)}{\beta!} \lambda^{3-2l} (-s)^l (-y)^\beta \mathcal{I}_5 (s,y) d\lambda dyds\\
	&-
	\int_{t/2}^t \int_{\mathbb{R}^2} \int_0^1
		(y \cdot \nabla) R^\bot R G(t-s,x-\lambda y)
		\mathcal{I}_5 (s,y)
	d\lambda dyds\\
	&-
	\sum_{2l+|\beta|=1}^3 \frac{\partial_t^l \nabla^\beta R^\bot RG(t,x)}{\beta!}
	\int_{t/2}^t \int_{\mathbb{R}^2}
		(-s)^l (-y)^\beta \mathcal{I}_5 (s,y)
	dyds.
\end{split}
\]
Here we used that $\int_{\mathbb{R}^2} \mathcal{I}_5 (s,y) dy = 0$ for the third and last parts.
%
Hence Hausdorff-Young inequality with Lemmas \ref{lem-sc} and \ref{RzG} and \eqref{scI} yields that
\[
\begin{split}
	\| J_3^\varepsilon (t) \|_{L^q (\mathbb{R}^2)}
	&\le
	C \int_\varepsilon^{t/2} (t-s)^{-\gamma_q-2} s^{-\frac12} ds
	+
	C t^{-\gamma_q - 2}
	\int_\varepsilon^{t/2} s^{-\frac12} ds
	+
	C \int_{t/2}^t (t-s)^{-\frac12} s^{-\gamma_q - 2} ds\\
	&
	+
	C \sum_{2l+|\beta| = 1}^3 t^{-\gamma_q-l-\frac{|\beta|}2} \int_{t/2}^t s^{-\frac52 + l + \frac{|\beta|}2} ds
\end{split}
\]
for $1 \le q \le \infty$.
Here the constants $C$ are independent of $\varepsilon$.
Thus, for any fixed $t > 0$, $\| J_3^\varepsilon (t) \|_{L^q (\mathbb{R}^2)}$ is bounded uniformly as $\varepsilon \to +0$.
Therefore Lebesgue convergence theorem concludes that $J_3$ is well-defined in $C((0,\infty),L^1 (\mathbb{R}^2) \cap L^\infty (\mathbb{R}^2))$.
Similar procedure with $\mathcal{I}_6$ instead of $\mathcal{I}_5$ guarantees the well-definedness of $J_4$.
The scaling properties are proved by elementary calculus, and \eqref{dcU} and \eqref{dcK} are already shown.
Then we complete the proof.
\end{proof}
\section{Proof of main results}
Far-field asymptotics of the first term on the right-hand side of \eqref{msu} are clarified by Lemma \ref{lem-asymplinsp}.
Hence, to prove Theorem \ref{thm-st}, we confirm \eqref{r4st} and \eqref{r4stinf}.
In Section 2, we showed Proposition \ref{prop-AsympVort} by using \eqref{decay-vort}.
In the similar way, we prove \eqref{r4st} and \eqref{r4stinf} by applying Proposition \ref{prop-AsympVort} instead of \eqref{decay-vort}.

\vspace{2mm}
\noindent
{\it Proof of Theorem \ref{thm-st}.}
The term of initial-data on \eqref{msu} is treated by Lemma \ref{lem-asymplinsp}.
We estimate for the nonlinear term.
Firstly, we consider the case that $\mu$ are large.
The error term $r_4$ on \eqref{expnl} is split to $r_4 (t) = r_4^1 (t) + r_4^2 (t) + \cdots + r_4^5 (t)$, where
\[
\begin{split}
	&r_4^j (t)\\
	&=
	\left\{\begin{split}
		\displaystyle
		&\iint_{Q_j}
			\biggl( R^\bot R G (t-s,x-y) - \sum_{l=0}^2 \frac{\partial_t^l R^\bot R G(t,x-y)}{l!} (-s)^l \biggr)
			\left( \omega u - \mathcal{I}_5 - \mathcal{I}_6 \right) (s,y)
		dyds,\\
		&\hspace{75mm}	j = 1,2,3,\\
		\displaystyle
		&\sum_{l=0}^2 \iint_{Q_j}
			\biggl( \frac{\partial_t^l R^\bot R G(t,x-y)}{l!}
			- \sum_{|\beta|=0}^{4-2l} \frac{\partial_t^l \nabla^\beta R^\bot R G(t,x)}{l! \beta!} (-y)^\beta \biggr)
			(-s)^l \left( \omega u - \mathcal{I}_5 - \mathcal{I}_6 \right) (s,y)
		dyds,\\
		&\hspace{80mm}	j = 4,5,
	\end{split}\right.
\end{split}
\]
and the separations $(0,t) \times \mathbb{R}^2 = Q_1 \cup Q_2 \cup Q_3 = Q_4 \cup Q_5$ are same as in \eqref{Q}.
By Taylor theorem, we see that
\[
	r_4^1 (t)
	=
	\int_0^{t/2} \int_{|y| > |x|/2} \int_0^1
		\frac{\partial_t^3 R^\bot R G(t-\lambda s, x-y)}{3!} \lambda^2 (-s)^3
		\left( \omega u - \mathcal{I}_5 - \mathcal{I}_6 \right) (s,y)
	d\lambda dyds.
\]
Thus, by Hausdorff-Young inequality with Lemmas \ref{lem-sc} and \ref{RzG} and Corollary \ref{cor-asympf}, we have that
\[
\begin{split}
	&\left\| |x|^\mu r_4^1 (t) \right\|_{L^1 (\mathbb{R}^2)}
	\le
	C \int_0^{t/2} \int_0^1 \| \partial_t^3 R^\bot R G(t-\lambda s) \|_{L^1 (\mathbb{R}^2)} \lambda^2 s^3
	\| |x|^\mu (\omega u - \mathcal{I}_5 - \mathcal{I}_6) (s) \|_{L^1 (\mathbb{R}^2)} d\lambda ds\\
	&\le
	C t^{- \frac52 + \frac\mu{2}} \log (2+t).
\end{split}
\]
On the same way, we obtain that $\| |x|^7 r_4^1 (t) \|_{L^\infty (\mathbb{R}^2)} \le C \log (2+t)$.
As long as $2 < \mu < 5$, the second part fulfills that
\[
\begin{split}
	&\| |x|^\mu r_4^2 (t) \|_{L^1 (\mathbb{R}^2)}
	\le
	C \int_0^t \int_0^1
		\| |x|^\mu \partial_t^3 R^\bot R G(t-\lambda s) \|_{L^1 (\mathbb{R}^2)}
		\lambda^2 s^3 \left\| (\omega u - \mathcal{I}_5 - \mathcal{I}_6) (s) \right\|_{L^1 (\mathbb{R}^2)}
	d\lambda ds\\
	&\le
	C \int_0^t \int_0^1 (t-\lambda s)^{-3 + \frac\mu{2}} \lambda^2 (1+s)^{-\frac12} \log (2+s) d\lambda ds
	\le
	C t^{-\frac52+\frac\mu{2}}\log (2+t)
\end{split}
\]
since $(t-\lambda s)^{-3+\frac\mu{2}}$ is integrable in $(\lambda,s) \in (0,1) \times (0,t)$.
We see that $\| |x|^7 r_4^2 (t) \|_{L^\infty (\mathbb{R}^2)} \le C \log (2+t)$ on the same way.
%
Similar procedure as above applied to the right-hand side of
\[
\begin{split}
	r_4^3 (t)
	&=
	\int_{t/2}^t \int_{|y| > |x|/2} \biggl( \int_0^1 \partial_t R^\bot R G(t-\lambda s,x-y) d\lambda (-s)\\
	&\hspace{15mm}
	- \sum_{l=1}^2 \frac{\partial_t^l R^\bot R G(t,x-y)}{l!} (-s)^l \biggr) \left( \omega u - \mathcal{I}_5 - \mathcal{I}_6 \right) (s)
	dyds
\end{split}
\]
provides that $\| |x|^\mu r_4^3 (t) \|_{L^1 (\mathbb{R}^2)} \le C t^{-\frac52 + \frac{\mu}2} \log (2+t)$ and $\| |x|^7 r_4^3 (t) \|_{L^\infty (\mathbb{R}^2)} \le C \log (2+t)$.
The fourth part is represented by
\[
	r_4^4 (t)
	=
	\sum_{2l + |\beta| = 5} \int_0^t \int_{|y| \le |x|/2} \int_0^1
		\frac{\partial_t^l \nabla^\beta R^\bot R G(t,x-\lambda y)}{l!\beta!} \lambda^{4-2l} (-s)^l (-y)^\beta
		\left( \omega u - \mathcal{I}_5 - \mathcal{I}_6 \right) (s,y)
	d\lambda dyds.
\]
Here we see from Lemma \ref{RzG} that $|x|^\mu \partial_t^l \nabla^\beta R^\bot R G \in L^1 (\mathbb{R}^2)$ for $0 \le \mu < 5$.
Hence
\[
\begin{split}
	&\| |x|^\mu r_4^4 (t) \|_{L^1 (\mathbb{R}^2)}
	\le
	C \sum_{2l + |\beta| = 5} \| |x|^\mu \partial_t^l \nabla^\beta R^\bot R G (t) \|_{L^1 (\mathbb{R}^2)}
	\int_0^t	
		s^l \| x^\beta ( \omega u - \mathcal{I}_5 - \mathcal{I}_6 ) (s) \|_{L^1 (\mathbb{R}^2)}
	ds\\
	&\le
	C t^{-\frac52+\frac\mu{2}} \left( \log (2+t) \right)^2.
\end{split}
\]
We confirm on the same way that $\| |x|^7 r_4^4 (t) \|_{L^\infty (\mathbb{R}^2)} \le C (\log (2+t))^2$.
Treatment for the last part is little complicated.
We estimate for $\| |x|^\mu r_4^5 (t) \|_{L^1 (\mathbb{R}^2)}$.
Since the mean value theorem provides that
\[
\begin{split}
	&r_4^5 (t)
	=
	-\int_0^t \int_{|y| > |x|/2}
		\biggl( \int_0^1 (y \cdot \nabla) R^\bot R G(t,x-\lambda y) d\lambda\\
		&\hspace{15mm} + \sum_{|\beta|=1}^4 \frac{\nabla^\beta R^\bot R G(t,x)}{\beta!} (-y)^\beta \biggr)
		\left( \omega u - \mathcal{I}_5 - \mathcal{I}_6 \right) (s,y)
	dyds\\
	&+ \sum_{l=1}^2 \int_0^t \int_{|y| > |x|/2}
		\biggl( \frac{\partial_t^l R^\bot R G(t,x-y)}{l!}
		 - \sum_{|\beta|=0}^{4-2l} \frac{\partial_t^l \nabla^\beta R^\bot R G(t,x)}{l!\beta!} (-y)^\beta \biggr)\\
		&\hspace{15mm} (-s)^l \left( \omega u - \mathcal{I}_5 - \mathcal{I}_6 \right) (s,y)
	dyds,
\end{split}
\]
we have that
\[
\begin{split}
	&\| |x|^\mu r_4^5 (t) \|_{L^1 (\mathbb{R}^2)}
	\le
	C \sum_{|\beta|=1}^4 \| |x|^{|\beta|-1} \nabla^\beta R^\bot R G(t) \|_{L^1 (\mathbb{R}^2)}
	\int_0^t	
		\| |x|^{\mu + 1} \left( \omega u - \mathcal{I}_5 - \mathcal{I}_6 \right) (s) \|_{L^1 (\mathbb{R}^2)}
	ds\\
	&+
	C \sum_{l=1}^2 \sum_{|\beta|=0}^{4-2l} \| |x|^{|\beta|} \partial_t^l \nabla^\beta R^\bot R G(t) \|_{L^1 (\mathbb{R}^2)}
	\int_0^t	
		s^l \left\| |x|^\mu \left( \omega u - \mathcal{I}_5 - \mathcal{I}_6 \right) (s) \right\|_{L^1 (\mathbb{R}^2)}
	ds\\
	&\le
	C t^{-\frac52 + \frac\mu{2}} (\log (2+t))^2.
\end{split}
\]
The estimate for $\| |x|^7 r_4^5 (t) \|_{L^\infty (\mathbb{R}^2)}$ does not require the mean value theorem, and we see that
\[
\begin{split}
	&\| |x|^7 r_4^5 (t) \|_{L^\infty (\mathbb{R}^2)}
	\le
	C \sum_{2l + |\beta| = 0}^{4} \| |x|^{|\beta|} \partial_t^l \nabla^\beta R^\bot R G (t) \|_{L^\infty (\mathbb{R}^2)}
	\int_0^t s^l \| |x|^7 (\omega u - \mathcal{I}_5 - \mathcal{I}_6) (s) \|_{L^1 (\mathbb{R}^2)} ds\\
	&\le
	C \log (2+t).
\end{split}
\]
Therefore we get \eqref{r4st} and \eqref{r4stinf} for large $\mu$.
Estimates for $\mu = 0$ are easier than the above.
Indeed, on the same way as in the proof of Proposition \ref{prop-defJ}, we see that
\[
\begin{split}
	&r_4 (t)
	=
	\int_0^{t/2} \int_{\mathbb{R}^2} \int_0^1
		\frac{\partial_t^3 R^\bot R G(t-\lambda s,x-y)}{3!} \lambda^2 (-s)^3
		\left( \omega u - \mathcal{I}_5 - \mathcal{I}_6 \right) (s,y)
	d\lambda dyds\\
	&+
	\sum_{2l + |\beta| = 5} \int_0^{t/2} \int_{\mathbb{R}^2} \int_0^1
		\frac{\partial_t^l \nabla^\beta R^\bot R G(t,x-\lambda y)}{l!\beta!} \lambda^{4-2l} (-s)^l (-y)^\beta
		\left( \omega u - \mathcal{I}_5 - \mathcal{I}_6 \right) (s,y)
	d\lambda dyds\\
	&-
	\int_{t/2}^t \int_{\mathbb{R}^2} \int_0^1
		(y \cdot \nabla) R^\bot R G (t-s,x-\lambda y)
		\left( \omega u - \mathcal{I}_5 - \mathcal{I}_6 \right) (s,y)
	d\lambda dyds\\
	&-
	\sum_{2l+|\beta|=1}^4 \frac{\partial_t^l \nabla^\beta R^\bot R G (t)}{l!\beta!} \int_{t/2}^t \int_{\mathbb{R}^2}
		(-s)^l (-y)^\beta \left( \omega u - \mathcal{I}_5 - \mathcal{I}_6 \right) (s,y)
	dyds.
\end{split}
\]
Hence, under the condition $|x|^5 \omega_0 \in L^1 (\mathbb{R}^2)$, Hausdorf-Young inequality with Lemmas \ref{lem-sc} and \ref{RzG} and Corollary \ref{cor-asympf} yields for $1 \le q \le \infty$ that
$
	\| r_4 (t) \|_{L^q (\mathbb{R}^2)}
	\le
	C t^{-\gamma_q - \frac52} ( \log (2+t) )^2.
$
Consequently we obtain \eqref{r4st} and \eqref{r4stinf} for $\mu =0$ and conclude the proof by H\"older embedding.\hfill$\square$

\vspace{3mm}

Also, by using a coupling of \eqref{decay} and \eqref{decay-vort} instead of Corollary \ref{cor-asympf} on the similar way, we see Proposition \ref{prop-lows}.
At last, we prove our second assertion.

\vspace{2mm}

\noindent
{\it Proof of Theorem \ref{thm-t}.}
The estimate for first term on \eqref{msu} is given by Lemma \ref{lem-asymplint}.
In Section 1, we expanded the nonlinear term as
\[
	\int_0^t R^\bot R G (t-s) * (\omega u) (s) ds
	=
	\sum_{m=1}^4 U_m^\infty (t) + \sum_{m=3}^4 (K_m + J_m + V_m) (t)
	+ (U_4^t - U_4^\infty) + r_4 + r_4^\infty.
\]
The last three parts are error terms.
The estimate for $r_4$ is given in the proof of Theorem \ref{thm-st} yet.
By employing Corollary \ref{cor-asympf}, we see \eqref{Usa} and \eqref{r4inft} for the third and last parts, respectively.
Here we used that $\| |x|^\mu (\mathcal{I}_p (1+t) - \mathcal{I}_p (t)) \|_{L^1 (\mathbb{R}^2)} \le C t^{-\frac{p}2 + \frac\mu{2}} (1+t)^{-1}$.
Thus we complete the proof.\hfill$\square$

\vspace{3mm}

Similar procedure as above are sure available in higher-dimensional cases.
However, it is little tough to introduce an asymptotic expansion since the vorticity is given by a tensor and some correction terms are required in those cases (cf.\cite{moto}).

Before closing this paper, we confirm that \eqref{decay-vort} is not strange.
The preceding works \cite{Kkvc, Kkvc-Trrs} provide the decay rates as $t \to +\infty$.
We treat the singularity as $t \to +0$.
For the first term on the right-hand side of \eqref{mcurl}, we see that $\| |x|^k G(t) * \omega_0 \|_{L^q (\mathbb{R}^2)} \le C (\| (|x|^k G(t)) * \omega_0 \|_{L^q (\mathbb{R}^2)} + \| G(t) * (|x|^k \omega_0) \|_{L^q (\mathbb{R}^2)})$.
The second term on the right-hand side yields the singularity and we see that $\| |x|^k G(t) *\omega_0 \|_{L^q (\mathbb{R}^2)} = O (t^{-\gamma_q})$ as $t \to +0$.
Hence the singularity of $\omega$ is coming from Gronwall type technique applied to \eqref{mcurl}.
If we use the condition that $|x|^{k+m} \omega_0 \in L^1 (\mathbb{R}^2)$ for some $m$ in the above, then the sigularity is mitigated.

\end{document}